\documentclass[12pt,a4paper]{article}
\usepackage[centertags]{amsmath}
\usepackage{amsfonts}
\usepackage{amssymb}
\usepackage{amsthm}
\usepackage{newlfont}
\usepackage{indentfirst}
\title{The Volume Entropy of a Riemannian Metric Evolving by the Ricci Flow on a Manifold of Dimension 3 or Above}
\author{Catalin C. Vasii}
\date{}
\theoremstyle{plain}
\newtheorem{te}{Theorem}[section]
\newtheorem{lm}[te]{Lemma}
\newtheorem{pr}[te]{Proposition}

\theoremstyle{definition}

\theoremstyle{remark}
\newtheorem{re}[te]{Remark}

\begin{document}
\maketitle
\typeout{:?0000}
 \begin{abstract}
In this paper it is proven that the volume entropy of a riemannian
metric evolving by the Ricci flow, if does not collapse,
nondecreases. Therefore it provides a sufficient condition
for a solution to collapse. Then, for the limit solutions of type I
or III, the limit entropy is the limit of the entropy as $t$
approaches the singular (finite or not) time.

\par{\bf key words and phrases:}  The Ricci Flow, Collapsing, Volume
entropy;
\par {\bf msc subj. class:} 58J35, 37B40

 \end{abstract}
 \section{Introduction}
 Consider on a compact manifold $M$ of dimension greater the 3 a
riemannian metric evolving by the Ricci Flow :
$\frac{\partial}{\partial t}g(t)=-2Rc_{g(t)}$. In \cite{manning} A.
Manning studies the volume entropy of a surface evolving under the
normalized Ricci Flow. Here, the author studies the volume entropy
of a manifold of dimension 3 or above, evolving by the unnormalized
flow, in connection with collapsing riemannian manifolds. In
\cite{per1} G. Perelman proves a noncollapsing theorem for a
solution to the Ricci flow that develops finite time singularities.
Despite this, a solution $g(t)$ that exists $\forall t$, could
collapse. The present paper deals with this case.

\par Let $\pi :\tilde M\to M$ be the universal cover of $M$.
Geometric data on $\tilde M$ will be denoted with tilde (such as
$\tilde R$ for the scalar curvature) to be distinguished from the
data on $M$. Consider The Ricci Flow on $\tilde M$:
$\frac{\partial}{\partial t}\tilde g(t)=-2\tilde{Rc}_{\tilde g(t)}$
and the one on $M$: $\frac{\partial}{\partial t}g(t)=-2Rc_{g(t)}$.
If $\tilde g(0)$ is the covering riemannian metric of $g(o)$ then,
as symmetries are preserved by the Ricci flow, $\tilde g(t)$ will be
the covering riemannian metric of $g(t)$, for all $t$ for which the
two solution exist.
\par By definition, the volume growth function of a compact
riemannin manifold $(M,g)$ is $\omega (M,g,r)=\frac{1}{r}\log Vol
\tilde B(r)$, where $\tilde B(r)$ is the radius $r$ ball in the
universal cover $(\tilde M, \tilde g)$ (the center of the ball
doesn't matter). The volume entropy is the quantity
$h(M,g)=lim_{r\to \infty}\omega (M,g,r)$.

 \section{Evolution of The Volume Entropy}
\begin{te}
If the injectivity radius $i(M,g(t))\geq i\in \mathbb R$ uniformly
in time then the volume entropy $h(M,g(t))$ is nondecreasing. If
$\int_{\tilde{B} (r(t))}\tilde{R} d\mu_{\tilde{g}(t)}-\int_{\tilde
S(r(t))}\tilde H d\mu_{\tilde \sigma(t)}\geq 0$ then $h(M,g(t))$ is again
nondecreasing.
\end{te}
\proof Consider the volume growth function of $(M,g(t))$, $\omega
(M,g(t),r)=\frac{1}{r}\log Vol \tilde B(r(t))$. We have that
$\frac{\partial}{\partial t}Vol \tilde B (r(t))=\lim_{h\to
0}\frac{1}{h}\bigl(\int _{\tilde B(r(t+h))} d\mu_{\tilde g(t+h)} -
\int _{\tilde B(r(t))} d\mu_{\tilde g(t)} \bigr)=\lim_{h\to
0}\frac{1}{h}\bigl(\int _{\tilde B(r(t+h))} d\mu_{\tilde g(t+h)} -
\int _{\tilde B(r(t+h))} d\mu_{\tilde g(t)} $ $+\int _{\tilde
B(r(t+h))} d\mu_{\tilde g(t)}$ $-\int _{\tilde B(r(t))} d\mu_{\tilde
g(t)} \bigr)=\lim_{h\to 0}\frac{1}{h}\bigl(\int _{\tilde B(r(t+h))}
d\mu_{\tilde g(t+h)} - \int _{\tilde B(r(t+h))} d\mu_{\tilde
g(t)}\bigr)+$ \linebreak $\lim_{h\to 0}\frac{1}{h}\bigl(\int
_{\tilde B(r(t+h))} d\mu_{\tilde g(t)} - \int _{\tilde B(r(t))}
d\mu_{\tilde g(t)}\bigr)=\lim_{h\to 0} \int _{\tilde
B(r(t+h))}\frac{1}{h}(d\mu_{\tilde g(t+h)}-d\mu_{\tilde g(t)})$
$+\lim_{h\to 0}\frac{1}{h}\int _{\tilde B(r(t+h))\backslash \tilde
B(r(t))}d\mu_{\tilde g(t)}$.
Since, as in \cite{ham1}, the
riemannian measure evolves  under the Ricci Flow by
\begin{equation}\label{volum}
\frac{\partial}{\partial t}d\mu_{\tilde g(t)}=-\tilde R d\mu_{\tilde
g(t)}
\end{equation}
we shall have $\frac{\partial}{\partial t}Vol \tilde B (r(t))=-\int
_{\tilde B(r(t))} \tilde R d\mu_{\tilde g(t)}+\lim_{h\to
0}\frac{\partial}{\partial h}\int _{\tilde B(r(t+h))\backslash
\tilde B(r(t))}d\mu_{\tilde g(t)}$ $=-\int _{\tilde B(r(t))} \tilde
R d\mu_{\tilde g(t)}+\frac{\partial}{\partial h}\lim_{h\to 0}\int
_{\tilde B(r(t+h))\backslash \tilde B(r(t))}d\mu_{\tilde g(t)}$.

Thus, since the second term vanishes, we obtain the evolution
\begin{equation}\label{evolutia}
\frac{\partial}{\partial t}\omega(M,
g(t),r)=-\frac{1}{r}\frac{1}{Vol(\tilde B(r(t)))} \int _{\tilde
B(r(t))} \tilde R d\mu_{\tilde g(t)}
\end{equation}
for the volume growth function. Now, the injectivity radius of the
universal riemannian cover satisfies $i (\tilde M,\tilde g(t))\geq
i(M,g(t))\geq i$ so, by a theorem from \cite{shen} that generalizes
 \cite{bg}  
we have that $$\lim_{r\to \infty}\frac{1}{Vol(\tilde B(r(t)))} \int _{\tilde
B(r(t))} \tilde R d\mu_{\tilde g(t)}\leq n(n-1)C$$ for some positive $C\in \mathbb R$, so
$$-\lim_{r\to \infty}\frac{1}{Vol(\tilde B(r(t)))} \int _{\tilde
B(r(t))} \tilde R d\mu_{\tilde g(t)}\geq -n(n-1)C.$$
Therefore, from
equation \ref{evolutia}, for $r \to \infty$, we get
$\frac{\partial}{\partial t}h(M, g(t))\geq 0$.
\par For the second assertion, we have that $\frac{\partial^2}{\partial r^2}\omega (M,g(t),r)=
\frac{2}{r^2}\omega(M,g(t),r)-\frac{1}{r^2 Vol \tilde
B(r(t))}-\frac{Vol^2 \tilde S (r(t))}{r Vol^2 \tilde B (r(t))}+
\frac{1}{r Vol \tilde B(r(t))}\int _{\tilde S(r(t))} \tilde H
d\sigma_{\tilde g(t)}$, where $d\sigma_{\tilde g(t)}$ is the
restriction  of the riemannian measure of the universal cover to the
sphere $\tilde S(r(t))$, and $\tilde H$ denotes the mean curvature
of the sphere. Thus, via equation (\ref{evolutia}), we get
\begin{eqnarray}\label{heateq}
\frac{\partial}{\partial
t}\omega(M,g(t),r)=\frac{\partial^2}{\partial r^2}\omega(M,g(t),r)-
\frac{2}{r^2}\omega(M,g(t),r)+ \frac{Vol \tilde S(r(t))}{r^2 Vol
\tilde B(r(t))}+ \cr
 \frac{Vol^2 \tilde S(r(t))}{r Vol \tilde B(r(t))}+ \frac{1}{r Vol \tilde B(r(t))}
\Bigl( \int_{\tilde{B} (r(t))}\tilde{R}
d\mu_{\tilde{g}(t)}-\int_{\tilde S(r(t))}\tilde H d\mu_{\tilde
\sigma(t)}\Bigr)
\end{eqnarray}
and this is a heat equation! Now, assuming by hypothesis that
$\int_{\tilde{B} (r(t))}\tilde{R} d\mu_{\tilde{g}(t)}-\int_{\tilde
S(r(t))}\tilde H d\mu_{\tilde \sigma(t)}\geq 0$ $\forall r$, in
equation (\ref{heateq}) above the last three terms are positive,
which implies that

\begin{equation}\label{supersol}
\frac{\partial}{\partial t}\omega(M,g(t),r)\geq
\frac{\partial^2}{\partial r^2}\omega(M,g(t),r)-
\frac{2}{r^2}\omega(M,g(t),r).
\end{equation}

In this point of the proof, we can apply the maximum principle for a
supersolution to the equation (\ref{supersol}). In order to do this,
consider the ODE $\frac{dx}{dt}=-\frac{2}{r^2}x$ which has the
solution $x=Ce^{(-2/r^2)t}$. By the maximum principle, the ODE gives
pointwise bounds to the PDE, so

$$\omega(M, g(t),r)\geq \omega(M, g(0),r)e^{(-2/r^2)t}$$
and this happens for all $r$. Taking limit as $r\to \infty$, we get

\begin{equation}\label{limita}
h(M,g(t))\geq h(M,g(0)).
\end{equation}
Let $t_k \nearrow \infty$. By repeating the argument above with
$t_k$ as origin, we obtain that $h(M, g(t))$ is nondecreasing.
 \hfill $\square$
\begin{re}
On a gradient Ricci soliton, the volume entropy is constant. Indeed,
a gradient Ricci soliton satisfies
$$-R= \Delta f +\frac{n\epsilon}{2t}$$
Where $\epsilon =-1,0,1$ according to the gradient Ricci soliton is
shrinking, steady or expanding. Then, by equation (\ref{evolutia}),
since by the divergence theorem we have $ \int _{\tilde B(r(t))}
\Delta f d\mu_{\tilde g(t)}=0$, follows that

$$\frac{\partial}{\partial t}\omega(M,
g(t),r)=\frac{1}{r}\frac{1}{Vol(\tilde B(r(t)))} \int _{\tilde
B(r(t))}  \Delta f +\frac{n\epsilon}{2t}d\mu_{\tilde
g(t)}=\frac{n\epsilon}{2tr}$$ thus, by taking limit as $r\to
\infty$, follows that $\frac{d}{dt}h(M,g(t))=0$.
\end{re}

 \section{The Volume Entropy and The Limit Solutions}
\begin{lm}\label{parabolic}
Under a parabolic rescaling of the solution to the Ricci flow the
volume growth function changes by
\begin{equation}\label{parabolic}
\omega(M, g_i(t),r)=\frac{n}{2}\frac{\log |Rm(x_i,t_i)|}{r}+ \omega
\bigl(M, g \bigl(t_i +\frac{t}{|Rm(x_i,t_i)|}\bigr), r\bigr).
\end{equation}
\end{lm}
\proof

Since a riemannian covering $\tilde M\to M$, beign a local isometry,
has $|\tilde Rm(y)|=|Rm(\pi(y))|$, we shall omit the tilde in the
following computation, to simplify notations.
\par Let $(x_i,t_i)$ a sequence of points $x_i\in M$ and times $t_i\nearrow
T\in(0,\infty]$ which converges to the singular time $T$. Then, a
parabolic rescaling of the metric is given by the formula
\begin{equation}
g_i(t):=|Rm(x_i,t_i)|g \bigl(t_i +\frac{t}{|Rm(x_i,t_i)|}\bigr).
\end{equation}
Hence, the volume of balls changes according to
$$\int_{B(x_i,r)}d\mu_{g_i(t)}=\int_{B(x_i,r)}|Rm(x_i,t_i)|^{n/2}d\mu_{g(t_i
+\frac{t}{|Rm(x_i,t_i)|})}$$ so

\begin{eqnarray}
\omega(M, g_i(t),r)=\frac{1}{r}\log Vol B_{g_i(t)}(x_i,r)= \cr
\frac{1}{r}\log (|Rm(x_i,t_i)|^{n/2}Vol B_{g(t_i
+\frac{t}{|Rm(x_i,t_i)|})}(x_i,r))
\end{eqnarray}

and equation (\ref{parabolic}) follows. \hfill $\square$
\par Let $(M_\infty, g_\infty (t))$ be the limit solution.
\begin{pr}
If $(M, g(t))$ is a type I or III solution to the Ricci flow, then
$h(M,g_\infty)=\lim_{t\to T} h(M, g(t))$.
\end{pr}
\proof Consider a type I solution, i.e. a solution that exhibits a
singularity in finite time $T$ and $(T-t)|Rm(x,t)|\lneq \infty$.
Assume $|Rm(x_i,t_i)|$ has not polynomial growth, i.e.
$|Rm(x_i,t_i)|\gneq \frac{1}{(T-t_i)^p}$ for some positive $p$,
hence $(T-t_i)^p
|Rm(x_i,t_i)|=(T-t_i)^{p-1}(T-t_i)|Rm(x_i,t_i)|\gneq 1$, so
$(T-t_i)^{p-1}\cdot \textrm{finite}\gneq 1$ and this contradicts the
fact that $t_i\to T$. Therefore $|Rm(x_i,t_i)|$ has polynomial
growth, and by equation (\ref{parabolic}) $$\lim_{i\to
\infty}\omega(M, g_i(t),r)=\lim_{i\to \infty}\omega \bigl(M,
g\bigl(t_i +\frac{t}{|Rm(x_i,t_i)|}\bigr), r\bigr)$$ so $h(M_\infty,
g_\infty (t))=\lim_{t\to T}h(M,g(t))$.

\par Consider now a type III solution, i.e. a solution that does not exhibit a
singularity in finite time, and $t|Rm(x,t)|\lneq \infty$. Assume
again $|Rm(x_i,t_i)|$ has not polynomial growth, i.e.
$|Rm(x_i,t_i)|\gneq t_i^p$ for some positive $p$ and repeat the
above reasoning. Follows that $h(M_\infty, g_\infty (t))=\lim_{t\to
\infty}h(M,g(t))$ which completes the proof. \hfill $\square$

 \begin{tabbing}
 C\u{a}t\u{a}lin C. Vasii\\
 Department of Mathematics, \\
 "Politehnica" University of  Timisoara,\\
 P-ta Victoriei 2, 300006 Timisoara, Romania,\\
 catalin@math.uvt.ro\\
 \end{tabbing}

\begin{thebibliography}{3}
\bibitem[GREEN]{bg} Green, L.W.,: {\it Aufwiedersehnsfl\"{a}chen},
Ann. Math {\bf 78} (1963), 289-299;
\bibitem[HAMILT1]{ham2} Hamilton, R.S.: {\it Three-Manifolds With Positive Ricci
Curvature}, J. Differential Geom. {\bf 17} (1982), no. 2, 255-306;

\bibitem[HAMILT2]{ham1} Hamilton, R.S.: {\it The Formation of Singularities in the Ricci
Flow},
 Surveys in Differential Geometry, Vol II (Cambridge MA, 1993)
 7-136;
\bibitem[MANN]{manning} Manning, A.: {\it The Volume Entropy of a surface decreases along the Ricci
Flow}, Ergod. Th. \& Dynam. Sys. {\bf 23} (2003), 1-6;
\bibitem[PER]{per1} Perelman, G.: {\it The Entropy Formula for The Ricci Flow and its Geometric
Applications} arXiv: math DG/0211159;
\bibitem[SHEN]{shen} Shen, Z.: {\it Conjugate radius and positive scalar curvature},
 Math. Z. {\bf 238} (2001), 431-439;

\end{thebibliography}
\end{document}